\documentclass{article}
\usepackage{a4}
\usepackage{amsmath}
\usepackage{amscd}
\usepackage{amssymb}
\usepackage[all]{xy}

\newcommand{\cD}{{\mathcal D}}

\newcommand{\cF}{\mathcal F}
\newcommand{\MA}[1]{{M_H{(#1)}}}

\newcommand{\M}{{M{(A)}}}

\newcommand{\dsum}{\displaystyle\sum}

\newcommand{\EA}{{\widehat{A}}}

\newcommand{\End}[1]{{\mathrm{End}}\:(#1)}

\newcommand{\Hom}[2]{{\mathrm{Hom}}\:(#1,#2)}
\newcommand{\EndX}[2]{{{\mathrm{End}}_{#1}\:(#2)}}
\newcommand{\HomX}[3]{{\mathrm{Hom}}_{#1}\:(#2,#3)}

\newcommand{\Ker}[1]{{\mathrm{Ker}}\:(#1)}

\newcommand{\Ann}[2]{{{\mathrm{Ann}_{#1}}\:({#2})}}

\newcommand{\lra}{\longrightarrow}
\newcommand{\ra}{\rightarrow}

\renewcommand{\Im}[1]{{\mathrm{Im}}\:(#1)}
\renewcommand{\lim}{{\displaystyle\varinjlim}}

\newcommand{\AH}{{A\# H}}
\newcommand{\ABt}{{A\#_\nu B}}

\newtheorem{theorem}{Theorem}[section]
\newtheorem{proposition}[theorem]{Proposition}
\newtheorem{lemma}[theorem]{Lemma}
\newtheorem{corollary}[theorem]{Corollary}
\newtheorem{definition}[theorem]{Definition}
\newtheorem{example}[theorem]{Example}

\newenvironment{proof}{{\sc Proof.}\ }{\ \rule{1ex}{1ex}\vspace{0.5truecm}}

\begin{document}

\title{A central closure construction for certain extensions. Applications to Hopf algebra actions.\thanks{Work supported by Funda\c{c}\~{a}o para a Ci\^{e}ncia e a Tecnologia through the Centro de Matem\'{a}tica da Universidade do Porto. Available as a PDF file from {http://www.fc.up.pt/cmup.}}}

\author{Christian Lomp}
\date{}

 \maketitle

\begin{abstract}
Algebra extensions $A\subseteq B$ where 
$A$ is a left $B$-module such that the $B$-action extends the multiplication in $A$
are ubiquitous. 
We encounter examples of such extensions in the study of group actions, group gradings or more general Hopf actions as well as 
in the study of the bimodule structure of an algebra.
In this paper we  are extending R.Wisbauer's method
of constructing the central closure of a 
semiprime algebra using its multiplication algebra
to those kinds of algebra extensions. 
More precisely if $A$ is a $k$-algebra and $B$ some subalgebra of $\End{A}$ that
contains the multiplication algebra of $A$, then the self-injective hull $\EA$ of 
$A$ as $B$-module becomes an $k$-algebra
provided $A$ does not contain any nilpotent $B$-stable ideals. 
We show that under certain assumptions $\EA$ can be identified with a subalgebra
of the Martindale quotient ring of $A$.
This construction is then applied to Hopf module algebras.
\end{abstract}

\section{Introduction}
Let $k$ be a commutative ring with unit. All $k$-algebras in this
paper are considered to be associative with unit. Unadorned tensor
products are taken over $k$ and $\End{-}$ resp. $\Hom{-}{-}$ refer
to $k$-linear maps. Let $A$ be a $k$-algebra. For any $a\in A$,
denote by $L_a$ the $k$-linear map $L_a \in \End{A}$ with
$L_a(x)=ax$ for all $x\in A$. Denote by $R_a\in \End{A}$ the
$k$-linear map $R_a(x)=xa$ for all $a\in A$. The $k$-subalgebra of
$\End{A}$ generated by the maps $L_a$ is denoted by $L(A)$. The
{\bf multiplication algebra} $\M$ of $A$ is the $k$-subalgebra of
$\End{A}$ generated by all maps $L_a$ and $R_a$; i.e.
$$ \M := \left\langle \{ L_a, R_a \mid a\in A\} \right\rangle\subseteq \End{A}.$$

The $k$-algebra $A$ is a cyclic left $\M$-module whose submodules
are precisely the two-sided ideals of $A$.

\begin{definition}\label{DefErweiterungModul} We say that an extension $A\subseteq B$ of $k$-algebras is an
{\bf extension with additional module structure $\varphi$}, if there exists
a ring homomorphism
$$ \varphi : B \rightarrow \End{A}$$
such that $\varphi(a) = L_a$ for all $a\in A$. 
We denote the
left $B$-module action on $A$ by $\cdot$, i.e. $$b\cdot a :=
\varphi(b)(a)$$ for all $a\in A$ and $b\in B$.
\end{definition}

Obviously $A$ becomes a cyclic left $B$-module. We will call a left
ideal $B$-stable if it is a $B$-submodule of $A$.
Let us denote by $\alpha$ the left $B$-linear map $\alpha: B
\longrightarrow A$ mapping an element $b$ of $B$ to $b\cdot 1_A$.

\begin{example}
Let $B:=A^e:=A\otimes A^{op}$ be the enveloping algebra of $A$ and
define a ring homomorphism $\varphi$ from $A^e$ to $\End{A}$ by
$\varphi(a\otimes b):= L_a\circ R_b$. Identifying $A$ with $A\otimes
1 \subseteq A^e$, we get that $A\subseteq A^e$ is an extension  with 
additional module structure $\varphi$. Note that $\Im{\varphi} = \M$.
\end{example}

\begin{example}\label{groupringexample}
Let $G$ be a group acting as ($k$-linear) automorphisms on $A$,
i.e. there exists a group homomorphism $\eta: G \rightarrow
Aut_k(A)$. We will use the notation $a^g:=\eta(g)(a)$ for all
$a\in A$ and $g\in G$. Define the skew group ring $A\# G$ whose
underlying $A$-submodule is the free left $A$-module with basis
$\{\overline{g}\mid g\in G\}$ and whose multiplication is given by
$(a\# \overline{g})(b\# \overline{h}) = a b^g \overline{gh}$. We might consider $A$ as a subring of
$A\# G$ by the map $A \rightarrow A\# G$ sending $a \mapsto a\# \overline{e}$
for all $a\in A$ where $e$ is the neutral element of $G$. An
action of $A\# G$ on $A$ is given by the ring homomorphism
$\varphi: A\# G\rightarrow \End{A}$ with $\varphi(a\# \overline{g}):= L_a
\circ \eta(g)$. Then $A\subseteq A\# G$ is an extension with additional module structure $\varphi$.
\end{example}

\begin{example}\label{derivationexample}
Let $\delta \in Der_k(A)$ be a $k$-linear derivation of $A$.
Consider the ring of differential operators $B=A[X;\delta]$, i.e.
as an $A$-module $B$ is equal to $A[X]$ but the multiplication is
constrained by $Xa - aX = \delta(a)$. Define a ring
homomorphism $\varphi: A[X;\delta]  \rightarrow \End{A}$ by
$\varphi(aX^n):=L_a\circ \delta^n$. Then $A\subseteq A[X; \delta]$ is an extension
with additional module structure $\varphi$.
\end{example}

\begin{example}
Let $H$ be an $k$-Hopf algebra acting on $A$. Denote the
action of an element $h\in H$ on $A$ by $\lambda_h\in \End{A}$.
The smash product $A\# H$ of $A$ and $H$ is the $A$-module $A\otimes H$ with multiplication given by
$(a\# h)(b\# g) := \sum_{(h)} a(h_1\cdot b) \# h_2g$ where $\Delta(h)=\sum_{(h)} h_1\otimes h_2$ is the comultiplication of $h$.
Define $\varphi:\AH\rightarrow \End{A}$ by
$\varphi(a\#h):=L_a\circ \lambda_h$. Then $A\subseteq A\# H$ is an extension with additional module structure $\varphi$.

For a group $G$ we might choose $H=k[G]$ and recover example \ref{groupringexample}.
For the trivial Lie algebra $\mathfrak{g}=k$ and its enveloping algebra $H=U(\mathfrak{g})=k[X]$  we recover example \ref{derivationexample}.
\end{example}

 If $A\subseteq B$ is an extension with additional module structure then $A\cap \Ann{B}{A} = 0$. Hence
$$A\subseteq B/\Ann{B}{A} \simeq \Im{\varphi} \subseteq \End{A}$$
is again an extension with additional module structure. Thus we
might replace $B$ by its image in $\End{A}$ and reduce ourselves
to extensions of $A$ inside $\End{A}$; where we identify $A$ with
$L(A)$. 

\begin{example}
Let $C$ be a $k$-bialgebra and $A$ a right $C$-comodule algebra 
with comodule structure $\rho:A\rightarrow A \otimes C$. 
For any $f\in C^*$ define an action on $A$ by $f\cdot a := (1\otimes f)\rho(a)$, for any $a\in A$. 
If we write $\rho(a)=\sum_{(a)} a_0 \otimes a_1$ then $f\cdot a = \sum_{(a)} a_0 f(a_1)$.
This defines an action of $C^*$ on $A$, i.e. we get a ring homomorphism $\phi: C^* \rightarrow \End{A}$.
Let $B$ be the subalgebra of $\End{A}$ generated by $L(A)$ and $\Im{\phi}$ then $A\subseteq B$ is an extension with
additional module structure. All left ideals of $A$ which are right $C$-comodules are $B$-stable. On the other hand if 
$C$ is a free $k$-module, then the $B$-stable left ideals of $A$ are precisely the left ideals which are right $C$-comodules.

As an application one might consider $G$-graded algebras $A$ (where $G$ is a monoid) as $k[G]$-comodule algebras.
In order to study the $G$-graded left ideals of $A$ one studies the $B$-stable left ideals of $A$ where 
$B$ is the subalgebra of $\End{A}$ generated by $L(A)$ and $(k[G])^*$.
\end{example}

If we want to investigate two-sided ideals that are stable
under a given action we have to restrict to extensions $A\subseteq B$ with additional module structure $\varphi$ 
such that $M(A) \subseteq \Im{\varphi}$. In some cases this happens
automatically. For instance let $A$ be an $k$-algebra with
involution $*$. Let $B$ be the subalgebra of $\End{A}$ generated
by $A$ and $*$. Since for any $a\in A$:
$$ R_a = *\circ L_{a^*}\circ *$$
we get $\M\subseteq B$. This means (as it is well-known) that any
left ideal of $A$ which is stable under $*$ is a two-sided ideal.
Note that $B$ can be seen as the factor ring of the skew-group
ring $A^e \# G$ where $G=\{id, g\}$ is the cyclic group
of order two and $g \in Aut(A^e)$ is given by
$$ (a\otimes b)^g := b^* \otimes a^*.$$
In this case we have that $A \subseteq A^e \# G$ is an extension with additional module structure.

Let $A\subseteq B$ be an arbitrary extension with additional module
structure. Recall the $B$-linear map $\alpha : B \rightarrow A$ with $(b)\alpha := b\cdot 1_A$ for all $b\in B$.
Note that $b\cdot a = b\cdot (a\cdot 1_A) = (ba)\cdot 1_A = (ba)\alpha$ for any $a\in A$ and $b\in B$.

We define the submodule of $B$-invariants of a left $B$-module as follows:
\begin{definition}
For any left $B$-module $M$ we denote
$$M^B := \{ m\in M\mid \forall b\in B, a \in A \:\: : \:\: b\cdot m = (b)\alpha \: m \}.$$
\end{definition}

Note that for any  $m\in M^B$ and any  $a\in A$ and $b\in B$ we have
$$b\cdot (am) = (ba)\cdot m = (ba)\alpha \: m = (b\cdot a) \: m.$$
The converse holds as well, i.e. if $b\cdot (am) = (b\cdot a)\: m$ for all $b\in B$ and $a\in A$ then $m\in M^B$.

 We can easily determine some elementary properties of ring extensions with additional module structure.
\begin{lemma}\label{ErweiterungModulstruktur} Let $A\subseteq B$ be an extension
with additional module structure. The the following properties hold:
\begin{enumerate}
\item[(1)] $B=A1_B \oplus \Ker{\alpha}$ as left $A$-modules.
\item[(2)] for all $M\in B$-Mod we have
$M^B =\Ann{M}{\Ker{\alpha}}$.
\item[(3)]
$\Psi:\EndX{B}{A}\longrightarrow A$ with  $\Psi(f):=(1_A)f$ is an
injective ring homomorphism with image $\Im{\Psi}=A^B$.
\item[(4)]
$\Psi_M : \HomX{B}{A}{M} \longrightarrow M^B$ with
$\Psi_M(f):=(1_A)f$ is an isomorphism of left $A^B$-modules.
\item[(5)] The isomorphisms $\Psi_M$ are natural transformations between the functors $\HomX{B}{A}{-}$ and $(-)^B$.
\item[(6)] The following statements are equivalent:
\begin{enumerate}
    \item[(a)] $()^B: B$-Mod $\lra A^B$-Mod is an exact functor.
    \item[(b)] $A$ is a projective left $B$-module.
    \item[(c)] there exists an element $t\in B^B$ with $(t)\alpha=1_A$.
    \item[(d)] there exists an idempotent $e\in B$ with $Be\simeq A$ as left $B$-modules and $eBe \simeq A^B$ as rings.
\end{enumerate}
\end{enumerate}
\end{lemma}

\begin{proof}
(1) Since $a'\cdot a=a'a$ for all $a,a' \in A$ holds, the map
$\beta: A \rightarrow B$ with $(a)\beta:=a1_B$ is $A$-linear
and lets $\alpha$ split as $A$-module homomorphism.
Hence as $A$-modules we have $B=A1_B \oplus \Ker{\alpha}$.\\
(2) If $m\in M^B$, then for all $b\in \Ker{\alpha}:$
$$b m = b(1_Am)=(b\cdot 1_A)m = (b)\alpha m = 0,$$
hence $M^B \subseteq \Ann{M}{\Ker{\alpha}}$.
On the other hand if $m \in \Ann{M}{\Ker{\alpha}}$, then
 $m\in M^B$, since from $B = A1_B \oplus \Ker{\alpha}$ it follows:
$$\forall b\in B, a\in A : b(am) = (ba)m = (ba)\alpha m = [b\cdot (a)\alpha]m = (b\cdot a)m,$$
(3 + 4) Let $f,g \in \EndX{B}{A}$, then $f$ and $g$
are in particular left $A$-linear and we have
$$\Psi(f\circ g):=(1_A)(f\circ g) := ((1_A)f)g = (1_A)f(1_A)g=\Psi(f)\Psi(g).$$
Thus $\Psi$ is a homomorphism of rings. Moreover for all
$b\in \Ker{\alpha}:$ $$b\cdot \Psi(f) = (b\cdot 1_A)f  = (b)\alpha
\circ f = 0.$$ By $(2)$ it follows that $\Im{\Psi}\subseteq A^B$.
On the other hand, for any $x\in A^B$, the right multiplication
$R_x$ is $B$-linear. To see this take any $a\in A$ and $b\in B$, then
$$(b\cdot a)R_x = (b\cdot a)x = (ba)\alpha x = (ba)\cdot x = b\cdot (ax) = b\cdot (a)R_x.$$
Hence $R_x \in \EndX{B}{A}$. Thus $A^B$ becomes a subring of  $A$
and every left $B$-module $M$ is also a left $A^B$-module with
$A^B$-submodule $M^B$. It follows as above that  $\Psi_M$ is an isomorphism.\\
(5) If $g:M\rightarrow N$ is a homomorphism between $B$-modules then
for any $f\in \HomX{B}{A}{M}$  we have 
$$(f)\Psi_M \circ g_{\mid M^B} = (1_A)f\circ g = (f\circ g)\Psi_N = (f)\HomX{B}{A}{g} \circ \Psi_N.$$
Thus $\Psi_M\circ g_{\mid M^B} = \HomX{B}{A}{g}\circ \Psi_N$, i.e. the isomorphisms $\Psi_M$ are natural transformations between the functors $(-)^B$ and $\HomX{B}{A}{-}$.\\
(6) $(a)\Leftrightarrow (b)$ holds by (5).\\
$(b)\Rightarrow (c)$ if $_BA$ is projective, then $\alpha$ splits and there exists a
 $B$-linear map $\beta: A\longrightarrow B$ with $\beta\alpha=id_A$. Set $t:=(1_A)\beta$. Then $t\in B^B$ by $(4)$ and 
 $(t)\alpha=1_A$.\\
$(c)\Rightarrow (b)$ if there exists an element $t\in B^B$ with $(t)\alpha = 1_A$,
then one defines $\beta: A \longrightarrow B$ as $(a)\beta =
at$. Since $t$ is in $B^B$,
 $\beta$ is $B$-linear and lets $\alpha$ split, i.e. $_BA$ is projective.
$(b)\Leftrightarrow (d)$ is clear.
\end{proof}

For $B=A^e$ we have $M^B = Z(M):=\{ m\in M\mid \forall a\in A \: am=ma \}$ and
$Z(-)$ is exact if and only if $A$ is a separable $k$-algebra.

For $B=A\# G$ we have $M^B = M^G := \{m\in M \mid \forall g\in G \: m^g = m \}$ and
$()^G$ is exact if and only if $G$ is finite and $A$ contains an element of trace one (property Lemma \ref{ErweiterungModulstruktur}.6(d)).

For $B=A[X;\delta]$ we have $M^B = \Ann{M}{X} = \{m\in M\mid Xm=0\}$.
It is impossible for $A$ to be a projective left $A[X;\delta]$-module simply because $B^B=l.ann_{A[X;\delta]}(X)=0$.

For $B=\AH$ we have $M^B=M^H:=\{m\in M\mid \forall h\in H \: h\cdot m = \epsilon(h)m \}$.
If $H$ is a finite dimensional Hopf algebra over a field $k$
then $()^H$ is exact in $\AH$-Mod if and only if there exists a left integral $t\in \int_l$ and an element $a\in A$ such that
$t\cdot a = 1$.

If $A$ is an algebra with involution $*$ and $B$ is the subalgebra generated by $A$ and $*$ in $\End{A}$, then
$M^B = Z(M;*):= \{m\in Z(M) \mid m^* = m \}$. Moreover $Z(-;*)$ is exact if and only if $A$ admits a separable idempotent
$\gamma = \sum_{i=1}^n x_i \otimes y_i$ such that $\gamma = \sum_{i=1}^n y_i^* \otimes x_i^*$.

\section{Prime and semiprime $B$-stable ideals}

For the rest of the paper we assume that 
$A\subseteq B$ is an extension with additional module structure $\varphi$
such that $\M\subseteq \Im{\varphi}$ (compare with \cite{Cabrera}).
Note that then every
$B$-submodule of $A$ is already a two-sided ideal. Moreover the
$B$-invariant elements $M^B$  for a left $B$-module $M$ are
$A$-centralizing. In  particular  $M^B \subseteq Z(M)$.

\begin{definition} A $B$-stable ideal $I$ of $A$ is called $B$-prime (resp. $B$-semiprime) if $KL\subseteq I$
(resp. $K^2\subseteq I$) implies $
K\subseteq I $ or $L\subseteq I$ (resp. $K\subseteq I$) for all $B$-stable ideals $K$ and $L$ of $A$.
$A$ is
called $B$-prime (resp. $B$-semiprime) if $0$ is a $B$-prime (resp.
$B$-semiprime) $B$-stable ideal. \end{definition}

If $I$ is a $B$-stable ideal of $A$, then there exists a ring homomorphism
$\varphi': B \longrightarrow \End{A/I}$ with
$M(A/I)\subseteq \Im{\varphi'}$. Let $B/I:=\varphi'(B)$, then $A/I
\subseteq M(A/I) \subseteq B/I$ is an extension with additional module structure $\varphi'$. With
this notation we prove easily the following Lemma:
\begin{lemma}
Let $A$ and $B$ as above. Let $P$ be a $B$-stable ideal of $A$ . Then $P$ is prime (resp. semiprime) if and only if
 $A/P$ is
$B/P$-prime (resp. $B/P$-semiprime).
\end{lemma}

Like in the classical case we have a description of $B$-stable semiprime ideals:

\begin{proposition}\label{SubdirektesProduktDarstellung} A $B$-stable ideal of $A$ is $B$-semiprime if and only if it is
the intersection of $B$-prime $B$-stable ideals.
\end{proposition}

\begin{proof}
$\Rightarrow$: Without loss of generality we might assume that $A$ is $B$-semiprime.
Let
$$I:=\bigcap \left\{ P \mid P \mbox{ is a prime $B$-stable ideal of $A$ } \right\}.$$
Assume $I \neq 0$. Then there exists $0\neq x_1 \in I$. Let
$I_1:= B\cdot x_1$ then $0\neq I_1 \subseteq I$. Since $A$ is $B$-semiprime, we have $(I_1)^2 \neq 0$. Hence
 $(I_1)^2$ contains a non-zero element $x_2$. Set $I_2 := B\cdot
x_2$. Again $(I_2)^2 \neq 0$, i.e. we may choose a non-zero element
$x_3 \in (I_2)^2$. Continuing this process we obtain a family $x_1,
x_2, x_3, \ldots $  of non-zero elements and a descending chain of non-zero $B$-stable ideals:
$$ I\supseteq I_1 \supseteq (I_1)^2 \supseteq I_2 \supseteq (I_2)^2
\supseteq \cdots \supseteq (I_{m-1})^2 \supseteq I_m \supseteq
\ldots$$ Consider the following set of $B$-stable ideals:
$$ {\mathcal Z}:= \{ P \subseteq A \mid P \mbox{ is a $B$-stable ideal of $A$ and for all } m: I_m \not\subseteq P \}.$$
Note that ${\mathcal Z}$ is not empty since $0 \in {\mathcal Z}$.
Let $\{ P_\lambda\}_\Lambda$ be an ascending chain of $B$-stable ideals in ${\mathcal Z}$.
Suppose  $I_m \subseteq \bigcup_{\Lambda} P_\lambda$ for some $m\geq 1$, then by definition of
$I_m = B\cdot x_m$ we have $x_m \in \bigcup_\Lambda P_\lambda$.
Thus there exists $\mu\in\Lambda $ such that $x_m \in P_\mu$ for some $\mu \in \Lambda$ and hence $I_m
\subseteq P_\mu$ - a contradiction to $P_\mu \in {\mathcal Z}$. Hence
$\bigcup_\Lambda P_\lambda\in {\mathcal Z}$ and we can apply Zorn's Lemma.
Let $P$ be a maximal element of ${\mathcal Z}$. We will show that  $P$ is a $B$-prime
$B$-stable ideal. Suppose there are  $B$-stable ideals $K,L$ such that  $KL \subseteq P$. Without loss of generality we might
assume $P\subseteq K$ and $P\subseteq L$.
If $L\neq P\neq K$, then by the maximality of  $P$ in ${\mathcal Z}$:
$K,L \not\in {\mathcal Z}$, i.e. there are
$m, n \geq 1$ with $I_m \subseteq K$ and $I_n \subseteq L$. Without loss of generality let
 $n \leq m$, then $$I_{m+1} \subseteq
(I_m)^2 \subseteq I_mI_n \subseteq KL \subseteq P,$$
a contradiction to $P\in {\mathcal Z}$. Hence $P=L$ or
$P=K$, i.e. $P$ is a  prime $B$-stable ideal. But this implies
$I\subseteq P$ and in particular $I_m \subseteq P$ for all $m$ - a contradiction.
Thus the intersection $I$ of all prime $B$-stable ideals is equal to zero.\\
The converse is clear: if $I^2=0$ for some $B$-stable ideal $I$ then $I^2 \subseteq P$ for any prime $B$-stable ideal $P$ of $A$.
Thus $I\subseteq P$ and hence $I\subseteq \bigcap P = 0$.
\end{proof}

We conclude that $A$ is $B$-semiprime if and only if $A$ is subdirect product of
$B/I$-prime algebras $A/I$.

Recall the following module theoretic notions:
The self-injective hull $\widehat{M}$ of a module $M$ is the largest $M$-generated submodule of the injective hull $E(M)$ of $M$, i.e.
if $E(M)$ denotes the injective hull of $M$, then 
$$\widehat{M} = \sum_{f\in \Hom{M}{E(M)}} \Im{f} = M\Hom{M}{\widehat{M}}.$$
The full subcategory of $R$-Mod consisting of submodules of $M$-generated modules is denoted by $\sigma[M]$.
The Lambek torsion theory in $\sigma[M]$ is the torsion category whose torsion class consists of all modules $X$ such that 
$\Hom{X}{\widehat{M}}=0$.
A submodule  $N$ of a module $M$ is
called {\bf dense}  if $M/N$ is a torsion module with respect to Lambek torsion theory, i.e.
$\Hom{M/N}{\widehat{M}}=0$. It is well-known that $N$ is dense in $M$ if and only if $\Hom{L/N}{M}=0$ for all submodules 
$N\subset L\subset M$ (see \cite[10.8]{wisbauer96}). $M$ is called  {\bf polyform}, if
every essential submodule of $M$ is $dense$. $M$ is called {\bf monoform}, if every non-zero submodule of $M$ is dense.
Dense submodules are also sometimes called rational (see \cite[chapter 10]{wisbauer96}).

\begin{lemma} \label{wesentlicheIdeale}
Let $A$ be $B$-semiprime. Then the following statements are
equivalent for a $B$-stable ideal $I$ of $A$:
\begin{enumerate}
\item[(a)] $I$ is a dense $B$-submodule of $A$.
\item[(b)] $I$ is an essential $B$-submodule of $A$.
\item[(c)] $JI\neq 0 \neq IJ$ for any non-zero $B$-stable ideal $J$ of $A$.
\end{enumerate}
\end{lemma}

\begin{proof}
$(a) \Rightarrow (b)$ dense submodules are essential (see \cite[chapter 10]{wisbauer96});\\
$(b)\Rightarrow (a)$ Let $K$ be a $B$-stable ideal of $A$ that contains $I$
and $f \in \HomX{B}{K/I}{A}$. Then $f(K/I)$ is a $B$-stable ideal of $A$.
Thus $N:=f(K/I)\cap I$ is a nilpotent $B$-stable ideal of $A$, since
$N^2\subseteq f(K/I)I=f(KI/I)=0$. Hence $N=0$, as $A$ is $B$-semiprime and $f(K/I)=0$,
as $I$ is an essential $B$-submodule. This shows $\HomX{B}{K/I}{A} = 0$, i.e.
$I$ is dense  in $A$.\\
$(b) \Leftrightarrow (c)$ For all
$B$-stable ideals $J$ we have: $(J\cap I)^2 \subseteq JI  \subseteq
J\cap I$. Since $A$ is $B$-semiprime we have $JI=0$ if and only if $I\cap J = 0$. Hence $I$ is an essential $B$-stable ideal
if and only if the left annihilator of $I$ does not contain any non-zero $B$-stable ideal. Analogously one concludes the same
statement for the right annihilator.
\end{proof}

As a corollary from the last Lemma we get:

\begin{corollary} \label{polyform} Let $A$ and $B$ be as above.
\begin{enumerate}
\item[(1)] If $A$ is $B$-semiprime, then  $A$ is a polyform
$B$-module and $A^B$ is reduced.
\item[(2)] If $A$ is $B$-prime,
then $A$ is a monoform $B$-module and $A^B$ is an integral domain.
\item[(3)] $A$ is $B$-prime if and only if $A$ is $B$-semiprime and a uniform $B$-module.
\end{enumerate}
\end{corollary}

\begin{proof}
(1) It follows from Lemma \ref{wesentlicheIdeale}[$(a)\Leftrightarrow (b)$]
that $A$ is polyform. Let $x \in A^B$ be such that $x^2=0$. Then $(Ax)^2=Ax^2=0$ shows that $Ax$ is a nilpotent $B$-stable ideal.
Thus $x=0$, i.e. $A^B$ is reduced. \\
(2) Let $A$ be $B$-prime and left $I$ be a non-zero $B$-stable ideal of $A$. Note that $JI=0$ implies $J=0$ for all $B$-stable ideals
as $A$ is $B$-prime and $I\neq 0$. By  Lemma \ref{wesentlicheIdeale} $I$ is a dense left $B$-submodule of $A$, i.e.
$A$ is monoform. If  $xy=0$ holds for $x,y \in A^B$, then $(Ax)(Ay)=Axy=0$. Since $A$ is $B$-prime, $x=0$ or $y=0$, i.e. $A^B$
is an integral domain.\\
(3) Follows from the definitions and Lemma \ref{wesentlicheIdeale}.
\end{proof}

\section{The central closure}

In the sequel we will extend Wisbauer's construction of the extended centroid and of the central closure of a semiprime algebra
(see \cite[chapter 32]{wisbauer96})
to our situation of an extension $A\subseteq B$ with additional module structure. 
We will reduce ourselves to subalgebras $B$ of $\End{A}$ which contain the multiplication algebra $\M$.

Let $Q_{max}(A)$ denote the maximal quotient ring of $A$.

\begin{theorem}\label{ErweitertesZentrum} Let $A$ be $B$-semiprime and let
$\EA$ be the self-injective hull of $A$ as $B$-module. Then the following hold:
\begin{enumerate}
\item[(1)] The map $\Psi: \EndX{B}{\EA} \longrightarrow
\EA^B$ with $\Psi(f):=(1_A)f$ is an isomorphism of $A^B$-modules and defines a ring structure on $\EA^B$ making it a
commutative, self-injective and von Neumann regular ring with subring $A^B$.
\item[(2)] There is a  bijection between the set of (essentially) closed $B$-stable ideals of $A$ and of the
set of central idempotents  of $\EA^B$.
\item[(3)] $\EA^B$ is a field if and only if $A$
is $B$-prime.
\item[(4)] $\EA^B$ is a finite product of $n$ fields  ($n\geq 1$) if and only if $A$ has finite Goldie dimension $n$ as left
$B$-module.
\item[(5)] If $A^B$ is large in $A$, i.e. $A^B\cap I \neq 0$ for all non-zero $B$-stable ideals $I$ of $A$
then
 $\EA^B=Q_{max}(A^B)$ and $A$ is non-singular as $A^B$-module
\end{enumerate}
\end{theorem}

\begin{proof}
(1) We know from Corollary \ref{polyform}, that $A$ is a polyform
$B$-module. Hence
\[\begin{CD}\EndX{B}{\EA} = \HomX{B}{A}{\EA}@>\Psi_{\EA}>>\EA^B\end{CD}\]
In particular $f=0$ if and only if $(1_A)f=0$ for all $f\in \EndX{B}{\EA}$.
$\EA^B$ carries a ring structure induced by $\Psi_{\EA}$ where
$$(1)f(1)g = (1)(f\circ g)$$
for all $f,g \in \EndX{B}{\EA}$. Moreover let $I:=(A)f^{-1}\cap (A)g^{-1} \cap
A$. Then for all $x,y \in I$ we have
$$ (xy)(f\circ g) = (x(y)f)g =(x)g(y)f = ((x)gy)f = (xy)(g\circ f),$$
i.e. $f\circ g - g\circ f \in \HomX{B}{\EA/I^2}{\EA}$. As intersection of two essential $B$-submodules, $I$ is essential
and no non-zero $B$-stable ideal annihilates $I$ on the left (see
Lemma \ref{wesentlicheIdeale}). Thus $I^2$ is also an essential $B$-submodule of $A$. By
Lemma \ref{wesentlicheIdeale} is $I^2$  dense. And henceforth as $\EA$ is polyform, $f\circ g = g\circ f$, i.e.
$\EndX{B}{\EA}\simeq \EA^B$ is  commutative. As endomorphism ring of a self-injective polyform module,
$\EA^B$ is self-injective and von Neumann regular and contains   $A^B$ as subring (see \cite[11.2]{wisbauer96}).\\
(2) follows from \cite[12.7]{wisbauer96};\\
(3) and (4) follow from (2) and (1);\\
(5) By \cite[11.5(1)]{wisbauer96} $\EA^B=Q_{max}(A^B)$. Let
$a\in A$ and $I$ an essential ideal of  $A^B$ with $aI=0$. Set
$J:=(B\cdot a)^B = (B\cdot a)\cap A^B$, Then $JI=0$ and hence $(J\cap I)^2=0$.
As $A^B$ is reduced and $I$ is essential in $A^B$ we conclude $J=0$. But since $A^B$ is large in $A$ we can also conclude 
$(B\cdot a)=0$, i.e.  $a=0$. Thus  $A$ is a  non-singular $A^B$-module.
\end{proof}

In the next theorem we will see that the self-injective hull $\EA$ itself carries a ring structure.

\begin{theorem}\label{zentraler Abschluss} Let $A$ be $B$-semiprime and let
$\EA$ be the self-injective hull of $A$ as $B$-module. Then
\begin{enumerate}
    \item[(1)] $\EA = Q_\cD(A)$ is the torsion theoretic quotient module with respect to the Lambek torsion theory
    $\cD$ in $\sigma[_BA]$ and $$\EA^B \simeq \lim \{ \HomX{B}{I}{A}| I \mbox{ is an essential $B$-submodule of $A$}\}.$$
    \item[(2)] The map $\phi: A \otimes_{A^B} \EndX{B}{\EA} \rightarrow \EA$ with $\psi(a\otimes f):=(a)f$ is left $B$-linear.
    Its kernel is an ideal and thus we might define a ring structure on $\EA$ given by the following multiplication:
          $$\forall a,b \in A; s,t \in \EA^B: (as)\cdot(bt) := (ab)st, $$
          where $A$ is a subring of $\EA$.
\end{enumerate}
Let  $\widehat{B}:=<B,\EA^B> \subseteq \End{\EA}$. Then
$\EA \subseteq M(\EA) \subseteq \widehat{B}$ is again an extension with additional module structure and the following hold:
\begin{enumerate}
    \item[(3)] $\EA$ is $\widehat{B}$-semiprime and a self-injective $\widehat{B}$-module.
    \item[(4)] $\EA$ is a non-singular $\EA^B$-module.
    \item[(5)] $A$ is $B$-prime if and only if $\EA$ is $\widehat{B}$-prime.
\end{enumerate}
We call $\EA$ the {\bf central closure} of $A$ with respect to $B$ and $\EA^B$ the {\bf extended centre} of
$A$ with respect to $B$.
\end{theorem}

\begin{proof}
(1) From the fact that $A$ is a polyform $B$-module it follows from \cite[9.13]{wisbauer96} $\EA=Q_\cD(A)$.
From \cite[9.17]{wisbauer96} follows the description of the endomorphism ring $\EA^B$.\\
(2) By definition the self-injective hull of a module $M$ is $M$-generated, i.e. 
the map $\phi: A \otimes \HomX{B}{A}{\EA} \rightarrow \EA$ with
$\phi(a\otimes f):=(a)f$ is an epimorphism of left $B$-modules,
where $B$ acting just on the first component of the tensor product $A\otimes \HomX{B}{A}{\EA}$. Since $A$ is polyform we have
$\HomX{B}{A}{\EA}=\EndX{B}{\EA}$. Let $a_i\in A$ and $f_i\in \End{\EA}$ and assume that
$\gamma:=\sum_{i=1}^n a_i \otimes f_i$ is in the kernel of $\phi$, i.e.
$$ 0=\phi(\gamma)=\sum_{i=1}^n (a_i)f_i$$
For any $b\otimes g$ we have
$$\phi((b\otimes g)\gamma)
= \phi\left(\sum_{i=1}^n ba_i \otimes gf_i\right)
= \sum_{i=1}^n (ba_i)g\circ f_i = b\left( \sum_{i=1}^n (a_i)f_i\right)g=0.$$
Hence $(b\otimes g)\gamma\in \Ker{\phi}$. Moreover
$$\phi(\gamma (b\otimes g))
= \phi\left(\sum_{i=1}^n a_ib \otimes f_ig\right)
= \sum_{i=1}^n (a_ib)f_i\circ g = \left( \sum_{i=1}^n (a_i)f_i\right)g b=0.$$
Thus $\gamma (b\otimes g)\in \Ker{\phi}$ shows that $\Ker{\phi}$ is an ideal of $A\otimes \EndX{B}{\EA}$
as claimed. This allows us to define an associative ring structure on the
$B$-module $\EA$ that is also compatible with that $B$-action and contains $A$ as a subring.\\
(3) By definition $\EA$ is a left $B$-module, hence there exists a ring homomorphism
$\Theta:B \longrightarrow \End{\EA}$ that is injective as $\EA$ is a faithful $B$-module.
Without loss of generality we might identify $B$ with its image $\Theta(B)$.
Let $\widehat{B}:=<\Theta(B), \EA^B>$. By hypothesis $M(A)\subseteq B$ implies
$$M(\EA) = M(A)\EA^B \subseteq \widehat{B}$$ and $\EA \subseteq \widehat{B}$ is an extension with additional module structure.
Let $I$ be a $\widehat{B}$-stable ideal of $\EA$, with $I^2=0$. In particular $(I\cap A)^2=0$ holds.
But since $I$ is $B$-stable, $A$ is essential as $B$-submodule of $\EA$ and $A$ is $B$-semiprime we conclude $I=0$.
Thus $\EA$ is $\widehat{B}$-semiprime.\\
Every $\widehat{B}$-endomorphism of $\EA$ is also $B$-linear.
On the other hand  $\EndX{B}{\EA}\simeq \EA^B \simeq \EndX{\widehat{B}}{\EA}$ holds. As $\EA$ was self-injective as $B$-module,
it is also self-injective as $\widehat{B}$-module.\\
(4) follows from \cite[11.11(5)]{wisbauer96}.\\
(5) Let  $I, J$ be non-trivial  $\widehat{B}$-stable ideals in
$\EA$. As $B\subseteq \widehat{B}$  these ideals are also $B$-stable. Since $A$ is essential as $B$-submodule, $(I\cap A)$ and
$(J\cap A)$ are non-trivial $B$-stable ideals of $A$ and $(I\cap A)(J\cap A)$ is contained in $IJ$. Hence if $A$ is $B$-prime, then
$\EA$ is $\widehat{B}$-prime. \\
On the other hand if $\EA$ is $\widehat{B}$-prime and $IJ=0$ for some $B$-stable ideals $I$ and $J$ of $A$,
then $(I\EA^B)(J\EA^B) = IJ\EA^B = 0$, i.e. $I\EA^B=0$ or $J\EA^B=0$. And $A$ is $B$-prime.
\end{proof}

For $B=\M$ we recover Wisbauer's construction of the central closure of $A$ (see \cite{wisbauer84}).

\section{The Martindale quotient ring}

Let $\cF$ denote the set of ideals of $A$ with zero left and right
annihilator. The {\it right Martindale ring of quotients} of $A$ is
$$ Q(A):=\lim \{ \HomX{-A}{I}{A} \mid I\in \cF\}.$$
Alternatively one might construct $Q(A)$ as follows: define an
equivalence relation $\sim$ on $\bigcup_{I\in \cF} \HomX{-A}{I}{A}$ by
letting $f:I\lra A$ be equivalent to $g:J\lra A$ if there exists 
$K\in \cF$ such that $K\subseteq I\cap J$ and $f_{\mid K} =
g_{\mid K}$. Denote by $[f]$ the equivalence class of a map $f:I\rightarrow A$.
Note that the equivalence class of the zero map $A\rightarrow A$
contains all maps $f$ that vanish on some ideal in $\cF$. Addition
is defined by $[f]+[g]:=[f+g: I\cap J \lra A]$ while multiplication is
set to be $[f][g]:=[fg: JI \lra A]$ where $fg$ denotes the composition map $a \mapsto f(g(a))$.

In case we have an extension $A\subseteq B$ with additional module structure, we are going to construct a subring
of $Q(A)$ related to all $B$-stable ideals in $\cF$.
Let  $\cF_B$ be the set of $B$-stable ideals with zero left and right annihilator.
{\bf We assume from now on that $A$ is $B$-semiprime and that the left annihilator of an $B$-stable ideal is again $B$-stable.}
By Lemma \ref{wesentlicheIdeale} all ideals in $\cF_B$ are essential $B$-submodules of $A$.
Consider the following construction:
$$Q_0(A):= \lim \{\HomX{-A}{I}{A} \mid I\in \cF_B\}.$$
We will refer to the elements of $Q_0(A)$ as equivalence classes
in the above sense, i.e.
$$ Q_0(A) = \left( \bigcup_{I\in \cF_B} \HomX{-A}{I}{A}\right)/\sim$$
where $$f\sim g \Leftrightarrow f_{\mid I} = g_{\mid I} \mbox{ for some }I\in \cF_B.$$
With the operations $+$ and $\cdot$ as above $Q_0(A)$ becomes a $k$-algebra and a subring of $Q(A)$.

Before we show that the central closure $\EA$ can be identified with a subring of 
$Q_0(A)$ we show that $\EA^B$ lies in the centre of $Q_0(A)$.

\begin{proposition}\label{ZenterOfMartindale} Let $A$ be $B$-semiprime and denote by $\EA$ the central closure of $A$ with 
respect to $B$. Assume that for any essential $B$-stable ideal $I$ of $A$ the left annihilator $l.ann_A(I)$ and 
right annihilator $r.ann_A(I)$ are $B$-stable ideals. Define for any $f\in \EndX{B}{\EA}$ the ideal $I_f:= (A)f^{-1}\cap A$. Then 
the map $$\psi: \EndX{B}{\EA} \rightarrow Q_0(A) \text{ with } f \mapsto [f : I_f \rightarrow A]$$
is an injective homomorphism of $k$-algebras whose image lies in the centre of $Q_0(A)$ and consists of
all elements $[f:I\rightarrow A]$ where $f$ is left $B$-linear.
\end{proposition}

\begin{proof}
For each endomorphism $f\in \EndX{B}{\EA}$ define $I_f := f^{-1}(A)\cap A.$ Since pre-images of
essential submodules are essential, $I_f$ is an essential
$B$-submodule of $A$. By Lemma \ref{wesentlicheIdeale}(c) and the hypothesis $I_f$ has zero left and right annihilator, i.e. $I_f\in \cF_B$.
We will show that $\psi$ is a ring homomorphism. Let $f,g \in \EndX{B}{\EA}$.
Note that $I_fI_g \in \cF_B$ and $I_fI_g \subseteq I_{fg}$
since for all $x\in I_f, y \in I_g$ the following holds :
$$fg(xy)=f(xg(y))=f(x)g(y)\in A.$$
Thus
$$\psi(f)\psi(g) = [f:I_f\ra A][g:I_g \ra A]
= [fg: I_fI_g \ra A] = [fg:I_{fg}\ra A]=\psi(fg).$$ This shows that
$\psi$ is a ring homomorphism. \\
Assume $\psi(f)=0$ for some $f\in \EndX{B}{\EA}$.
Then there exists an $J\in \cF_B$ with $J\subseteq I_f$ and
$f(J)=0$. Hence $f\in \HomX{B}{\EA/J}{\EA}=0$ as $J$ is dense by
Lemma \ref{wesentlicheIdeale}. Thus $f=0$ and  $\psi$ is injective.\\
Let $[f:I\rightarrow A]$ be such that $f$ is $B$-linear, then $f$ can be uniquely extended to an endomorphism $\overline{f}\in \EndX{B}{\EA}$
since $\EA$ is self-injective and polyform as $B$-module. By definition
$\psi(\overline{f}) = [f:I\rightarrow A]$ since $I\subseteq I_{\overline{f}}$ and $f_{\mid I} = \overline{f}_{\mid I}$.
Hence the image of $\psi$ consists of all elements $[f:I\rightarrow A]$ such that $f$ is $B$-linear.
\end{proof}

Let $\imath :A \rightarrow Q_0(A)$ be the inclusion of $A$ into $Q_0(A)$
 given by $\imath(a):=[L_a:A\rightarrow A]$.
Together with $\psi$ we have a map $A\times \EndX{B}{\EA} \rightarrow Q_0(A)$ sending a pair $(a,f)$ to the product
$\imath(a)\psi(f)$. Since $A^B\simeq \imath(A^B)\subseteq Z(Q_0(A))$  this map is $A^B$-balanced and induces
an $k$-algebra homomorphism $$\psi^*: A\otimes_{A^B} \EndX{B}{\EA} \rightarrow Q_0(A).$$
Recall from Theorem \ref{zentraler Abschluss} the $k$-algebra homomorphism
$\phi:A \otimes_{A^B} \EA^B \rightarrow \EA$ with $\phi(a\otimes f)=af$.

\begin{lemma} For any element $\gamma  \in A\otimes_{A^B} \EndX{B}{\EA}$ there exists a $B$-stable ideal $I\in \cF_B$ such that
$$\psi^*(\gamma)=[L_{\phi(\gamma)} : I \rightarrow A]$$
where $L_{\phi(\gamma)}(x)=\phi(\gamma) x$ for all $x\in I$.
\end{lemma}
\begin{proof}
Write $\gamma = \sum_{i=1}^n a_i\otimes x_i$. By definition
\begin{eqnarray*} \psi^*(\gamma)&=&\sum_{i=1}^n [L_{a_i}:A\rightarrow A][L_{x_i}:I_{x_i}\rightarrow A]\\
&=& \sum_{i=1}^n [L_{a_ix_i}:I_{x_i}\rightarrow A] \\
&=& \left[ \sum_{i=1}^n L_{a_ix_i}: \bigcap_{i=1}^n I_{x_i}\rightarrow A \right] \\
&=& [L_{\sum_{i=1}^n a_ix_i}: I \rightarrow A]
\end{eqnarray*}
where we set $I:=\bigcap_{i=1}^n I_{x_i} \in \cF_B$. Since $\phi(\gamma)=\sum_{i=1}^n a_ix_i$ our claim is proved.
\end{proof}

In particular this implies that $\Ker{\phi}\subseteq \Ker{\psi^*}$, i.e.
$\psi^*$ extends to an $k$-algebra homomorphism $\overline{\psi}:\EA\rightarrow Q_0(A)$.

Let $\jmath:A\rightarrow \EA$ denote the inclusion map.

\begin{proposition} The following diagram in the category of $k$-algebras commutes:

\begin{center}\xymatrix{
A \ar[rr]^{\imath}\ar[dr]_{\jmath} & & Q_0(A)\\
& \EA \ar[ur]^{\overline{\psi}} &
}\end{center}

The image $\Im{\overline{\psi}}$ is the subalgebra of $Q_0(A)$ generated by
the image of $A$ and all elements $[f:I\rightarrow A]$ such that $f$ is $B$-linear.
The kernel of $\overline{\psi}$ is equal to
$$ \Ker{\overline{\psi}} = \phi(\Ker{\psi^*}) = \bigcup_{I\in \cF_B} l.ann_{\EA}(I).$$
\end{proposition}

\begin{proof}
Assume that $\gamma \in \Ker{\psi^*}$, then $\psi^\ast(\gamma)=[L_\phi(\gamma):I \rightarrow A] = 0$ for some $I\in \cF_B$,
i.e. $\phi(\gamma)I=0$. Thus
$\phi(\gamma)$ is an element of the left annihilator in $\EA$ of $I$. 
Hence  $\Ker{\psi} \subseteq  \bigcup_{I\in \cF_B} l.ann_{\EA}(I).$
On the other hand each element $z$ in $l.ann_{\EA}(I)$ is mapped by $\overline{\psi}$ to the zero class in $Q_0(A)$.
\end{proof}

Under some conditions $\overline{\psi}$ is injective.

\begin{theorem}\label{CentralClosureEmbedds}
Let $A$ be $B$-semiprime and $\EA$ its central closure with respect to $B$. Assume that the following conditions are fulfilled:
\begin{enumerate}
  \item[(i)]  $l.ann_A(I)$ and $r.ann_A(I)$ are $B$-stable for all essential $B$-stable ideals $I$ of $A$.
  \item[(ii)] $l.ann_{\EA}(I)$ is $B$-stable for all $B$-submodules $I$ of $\EA$.
\end{enumerate}
Then $\overline{\psi} : \EA \rightarrow Q_0(A)$ is an injective $k$-algebra homomorphism.
\end{theorem}

\begin{proof} It is enough to show that $\phi(\Ker{\psi})=0$ since 
then $\Ker{\psi}=\Ker{\phi}$ holds, i.e. $\overline{\psi}$ is injective.
Let $I\in \cF_B$. Then $0=l.ann_A(I) = A \cap l.ann_{\EA}(I)$. Since $I$ is an essential $B$-submodule of $A$ and hence of $\EA$,
by hypothesis $l.ann_{\EA}(I)$ is $B$-stable. But then $l.ann_{\EA}(I)=0$ as $A$ is essential as $B$-submodule of $\EA$.
\end{proof}

\section{Applications to Hopf actions}

In this section we are going to apply our construction to Hopf actions.
Let $H$ be a Hopf algebra over $k$ and let $A$ be a left $H$-module algebra.
Then there exists a left $H$-module structure on $A$ given by some ring homomorphism
$\varphi: H\rightarrow \End{A}$. Let us denote the action of an element $h\in H$ as an endomorphism of $A$ by
$\lambda_h$, i.e.
$\lambda_h(a)=h\cdot a$ for all $a\in A$.
In order for $A$ to be a left $H$-module algebra the $H$-action has to satisfy the following
condition in $\End{A}$ for $h\in H$ and $a\in A$:
$$ \lambda_h\circ L_a = \sum_{(h)} L_{h_1\cdot a}\circ \lambda_{h_2}$$
where $\Delta(h)=\sum_{(h)} h_1 \otimes h_2$ denotes the comultiplication of $h$ in Sweedler notation.

The $k$-subalgebra generated by the maps $L_a, R_a$ and $\lambda_h$  for $a\in A$ and $h\in H$ is denoted by $\MA{A}$, i.e.
$$\MA{A} := \left\langle \{L_a, R_a, \lambda_h \mid a\in A, h\in H\} \right\rangle \subseteq \End{A}.$$

Instead of $\MA{A}$-prime resp. $\MA{A}$-semiprime one says that $A$ is $H$-prime resp. $H$-semiprime.
It is easy to verify that $A^{\MA{A}}=Z(A)\cap A^H$ holds. Let us denote $Z(A)\cap A^H$ by $Z(A)^H$.
Note that in general $Z(A)$ is not closed under the action of $H$, but $Z(A)$
is a left $H$-module algebra in case $H$ is cocommutative.

Let us first realize $\MA{A}$ as the factor of some kind of smash product.
If $A$ is commutative, then $\MA{A}$ is generated by $\{L_a, \lambda_h | a\in A, h\in H\}$.
Hence we might identify $\MA{A}$ with the
image of $\AH \rightarrow \End{A}$ mapping $a\#h$ to $L_a\circ \lambda_h$. In this case
$\MA{A}$ is isomorphic to $\AH/\Ann{\AH}{A}$.

In the following theorem we represent $\MA{A}$ as a factor ring of some smash product.
For that reason we are going to introduce a general construction:

\begin{definition} Let $A$ and $B$ be $k$-algebras with multiplication maps
$\mu_A$ resp. $\mu_B$. Let $\nu: B\otimes A \longrightarrow
A\otimes B$ be an $k$-linear map and define:
$$\begin{CD}\mu : (A\otimes B)\otimes (A\otimes B)
@>{1\otimes \nu\otimes 1}>> (A\otimes A) \otimes (B \otimes B)
@>{\mu_A \otimes \mu_B}>> A\otimes B.\end{CD}$$ If $A\otimes B$ becomes through $\mu$
 an associative $k$-algebra with unit $1_A \otimes 1_B$, then we will write  $\ABt$ and call this ring the
{\bf smash-product} or factorization structure of $A$ and  $B$ with respect to $\nu$.
\end{definition}

Caenepeel, Ion, Militaru and Zhu gave a characterisation of smash-products:

\begin{theorem}[{\cite[Theorem 2.5]{CaenepeelIonMilitaruZhu}}]\label{CharakterisierungSmashProdukt}
Let $A$,$B$ and $\nu$ as above.
Then  $\nu$ defines a smash-product for $A$ and $B$ if and only if the following statements hold:
\begin{enumerate}
    \item[(i)] $\nu(b\otimes 1_A) = 1_A \otimes b$ and $\nu(1_B \otimes a)=a \otimes 1_B$ for all $a\in A , b\in B$.
    \item[(ii)] The following diagrams commute:
\end{enumerate}

\hspace{-3mm}\xymatrix{
B\otimes B \otimes A \ar[r]^{\mu_B\otimes 1}\ar[dd]_{1\otimes \nu} & B\otimes A \ar[dr]^{\nu}&&B\otimes A\ar[dl]_{\nu}&B\otimes A \otimes A\ar[l]_{1\otimes \mu_A}\ar[dd]^{\nu\otimes 1} \\
&&A\otimes B&&\\
B\otimes A \otimes B\ar[r]_{\nu\otimes 1} &A\otimes B \otimes B\ar[ur]_{1\otimes\mu_B}&&A\otimes A\otimes B\ar[ul]^{\mu_A\otimes 1}&A\otimes B\otimes A\ar[l]^{1\otimes \nu} \\
}
\end{theorem}

Let $A^e:=A\otimes A^{op}$ be the enveloping algebra of $A$.
We want to define a smash product of $A^e$ and $H$.
If $H$ is cocommutative, then $A^e$ is a left $H$-module algebra and we can use the ordinary smash product
$A^e \# H$, but in general $A^e$ will not be an $H$-module algebra.

Define the map $\nu: H\otimes A^e \longrightarrow A^e \otimes H$ by
$$\nu (h\otimes a\otimes b ) := \sum_{(h)} (h_1\cdot a)\otimes (h_3\cdot b) \otimes h_2$$
for all $a,b \in A$ and $h\in H$.

Then $A^e {\#_\nu} H$ is a smash product in the above sense.
To see this we have to check that the diagrams above commute.

 Property (i) of Theorem \ref{CharakterisierungSmashProdukt} is obviously fulfilled. Let $a,b,x,y \in A$ and $h,g\in H$. Then
\begin{eqnarray*}
(1\otimes \mu_H)(\nu\otimes 1)(1\otimes \nu)&&(h\otimes g \otimes
(a\otimes b))\\
&&=\sum_{(g)} (1\otimes \mu_H)(\nu\otimes 1)(h\otimes (g_1\cdot a \otimes g_3\cdot b) \otimes g_2)\\
&&=\sum_{g,h)} (1\otimes \mu_H)((h_1g_1\cdot a)\otimes (h_3g_3\cdot b) \otimes h_2 \otimes g_2) \\
&&=\sum_{g,h)} h_1g_1\cdot a\otimes h_3g_3\cdot b \otimes h_2g_2.\\
&&=\nu(\mu_H \otimes 1)(h\otimes g \otimes (a\otimes b))
\end{eqnarray*}
Hence $(1\otimes \mu_H)(\nu\otimes 1)(1\otimes \nu)=\nu(\mu_H
\otimes 1)$ holds, i.e. the left part of the diagram in
Theorem \ref{CharakterisierungSmashProdukt}(ii) commutes. 

\begin{eqnarray*}
 (\mu_{A^e}\otimes 1)(1 \otimes \nu)(\nu\otimes 1)&&(h\otimes (x\otimes y) \otimes (a\otimes b))\\
&&=\sum_{(h)}(\mu_{A^e}\otimes 1)(1 \otimes \nu)(( h_1\cdot x \otimes h_3\cdot y)\otimes h_2 \otimes (a\otimes b))\\
&&=\sum_{(h)}(\mu_{A^e}\otimes 1)((h_1\cdot x \otimes h_5\cdot y) \otimes (h_2\cdot a\otimes h_4\cdot b) \otimes h_3)\\
&&=\sum_{(h)}(h_1\cdot x)(h_2\cdot a) \otimes (h_4\cdot b)(h_5\cdot y) \otimes h_3\\
&&=\sum_{(h)}((h_1\cdot (xa)) \otimes (h_3\cdot (by)) \otimes h_2)\\
&&=\nu(1\otimes \mu_{A^e}) (h\otimes (x\otimes y) \otimes (a\otimes b))
\end{eqnarray*}

Hence $(\mu_{A^e}\otimes 1)(1 \otimes \nu)(\nu\otimes
1)=\nu(1\otimes \mu_{A^e})$, i.e. the right part of the diagram of
Theorem \ref{CharakterisierungSmashProdukt}(ii) commutes and
$A^e {\#_\nu}H$ is a smash product.

One could also define a smash product on the $k$-algebras $A^{op}$ and $\AH$.
It is not difficult to check that the map $\sigma: A^{op}\otimes \AH \rightarrow \AH \otimes A^{op}$ with
    $$\sigma (b\otimes a\# h ) := \sum_{(h)} a\# h_1 \otimes (S(h_2) \cdot b),$$
will define a smash product $(\AH) {\#_\sigma} A^{op}$.
Moreover one checks that the map
$\Psi: (\AH){\#_\sigma}A^{op} \longrightarrow  A^e{\#_\nu}H$ with
$\Psi(a\# h \otimes x):=\sum_{(h)}(a \otimes (h_2\cdot x)) \# h_1$ is an isomorphism of $k$-algebras.

We can now represent $\MA{A}$ as the factor of the smash product $A^e{\#_\nu}H$.

\begin{theorem}
The map $\Phi : A^e {\#_\nu} H \longrightarrow \MA{A}$
with $\Phi((a\otimes b) \# h):= L_a\circ R_b \circ \lambda_h$ is a surjective map of $k$-algebras.
Moreover $A\subseteq A^e {\#_\nu} H $ is a ring extension with additional module structure $\Phi$.
\end{theorem}

\begin{proof}
We have the well-defined maps $L_{A^e}:A^e \longrightarrow
\End{A}$, and $\lambda_H:H \longrightarrow \End{A}$. Let $\mu$
denote the multiplication in $\End{A}$, then $\Phi =
\mu\circ(L_{A^e} \otimes \lambda_H)$ is well-defined. By definition
$\Im{\Phi}=\MA{A}$, i.e. $\Phi$ is surjective .
To show that $\Phi$ is a ring homomorphism note that for all $h\in H, a,b,x \in A$
$$ h\cdot (axb)=\dsum_{(h)} (h_1\cdot a)(h_2\cdot x)(h_3\cdot b)$$
holds and therefore also
$$\lambda_h\circ L_a\circ R_b = \dsum_{(h)} L_{h_1 \cdot a}\circ R_{h_3 \cdot b} \circ \lambda_{h_2}.$$
By definition of the multiplication in $A^e \#_\nu H$ this implies that
$\Phi$ is a ring homomorphism.
\end{proof}

In case $H$ is cocommutative, we have that $A^e \#_\nu H$ and $A^e \# H$ coincide:

\begin{proposition}
If $H$ is cocommutative, then $A^e$ is a left $H$-module algebra and
    $A^e \#_\nu H$ is equal to the ordinary smash product $A^e \# H$ of a module algebra and the Hopf algebra.
\end{proposition}

\begin{proof}
$A^e$ is always a left $H$-module by the diagonal module structure, i.e.
$$h\cdot (a\otimes b):=\sum_{(h)} (h_1\cdot a)\otimes (h_2\cdot b)$$
for all $h\in H$ and $a,b\in A$.
Suppose that $H$
is cocommutative.
Let $a\otimes x, b\otimes y \in A^e$ and $h\in H$. Then
\begin{eqnarray*}
h\cdot ((a\otimes x)(b\otimes y)) &=& h\cdot (ab \otimes yx) \\
&=& \sum_{(h)} h_1\cdot (ab) \otimes h_2\cdot (yx) \\
&=& \sum_{(h)} (h_1\cdot a)(h_2\cdot b) \otimes (h_3\cdot y)(h_4\cdot x) \\
&=& \sum_{(h)} (h_1\cdot a)(h_3\cdot b) \otimes (h_4\cdot y)(h_2\cdot x) \\
&=& \sum_{(h)} [(h_1\cdot a)\otimes (h_2\cdot x)][(h_3\cdot b) \otimes (h_4\cdot y)] \\
&=& \sum_{(h)} [h_1 \cdot (a \otimes x)][h_2 \cdot (b \otimes y)]
\end{eqnarray*}
Moreover $h\cdot (1\otimes 1) = \varepsilon(h)(1\otimes 1)$ holds, i.e.
$A^e$ is a left $H$-module algebra.\\
Moreover for all $h\in H$ and $a,x \in A$:
$$\nu(h\otimes (a\otimes x)) = \sum_{(h)} (h_1\cdot a) \otimes (h_3\cdot x) \otimes  h_2
= \sum_{(h)} [h_1\cdot (a\otimes x)] \otimes h_2.$$ Thus
the multiplication defined by $\nu$ is equal to the multiplication in the ordinary smash product,
i.e. $A^e\#_\nu H = A^e\# H$.
\end{proof}

Now we are in position to apply our previous results to Hopf module algebras.

\begin{theorem}\label{AnwendungZentralerAbschluss} Let $H$ be a $k$-Hopf algebra, $A$ be
an $H$-semiprime left $H$-module algebra and let
$\EA$ be the self-injective hull of $A$ as $A^e {\#_\nu} H$-module.
Then the following statements hold:
\begin{enumerate}
    \item[(1)] $A$ is a polyform $A^e {\#_\nu} H$-module and a subdirect product of $H$-prime module algebras.
    Furthermore $Z(A)^H$ is reduced.
    \item[(2)] $\EA$ is an $H$-semiprime left $H$-module algebra with submodule algebra $A$.
        $\EA$ is self-injective as ${\EA}^e {\#_\nu} H$-module and
        a non-singular module over the self-injective and von Neumann regular ring  $Z(\EA)^H$.
    \item[(3)] If $Z(A)^H$ is large in $A$, then
        $Z(\EA)^H=Q_{max}(Z(A)^H)$ and $A$ is non-singular as $Z(A)^H$-module
\end{enumerate}
\end{theorem}

\begin{proof}
Note that the $A^e {\#_\nu} H$-module structure of $A$ coincides with that of $\MA{A}$ since
$\MA{A} \simeq A^e {\#_\nu} H / \Ann{A^e {\#_\nu} H}{A}$.\\
(1) follows from Corollary \ref{polyform} and Proposition \ref{SubdirektesProduktDarstellung}.\\
(2) follows from Theorem \ref{zentraler Abschluss} and
Theorem \ref{ErweitertesZentrum}. Note that
$$\MA{\EA} = <\MA{A}, Z(\EA)^H> = \widehat{\MA{A}} \subseteq \End{\EA}.$$
We still have to prove that $\EA$ is a left
$H$-module algebra. The $H$-module structure on 
$\EA=A Z(\EA)^H$ is given by
$h\cdot(as)=(h\cdot a)s$. Let $as, bt \in \EA$ and
$h \in H$, then: {\small $$h\cdot [(as)(bt)] = (h\cdot
(ab))st = \sum_{(h)} (h_1\cdot a)(h_2 \cdot b)st = \sum_{(h)}
[(h_1\cdot a)s][(h_2\cdot b)t] = \sum_{(h)} [h_1 \cdot (as)][h_2
\cdot (bt)].$$} (3) follows from Theorem \ref{ErweitertesZentrum}.
\end{proof}

From \ref{ErweitertesZentrum}(2,3,4) follows also:

\begin{corollary}
Let $H$ be a $k$-Hopf algebra, $A$ be
an $H$-semiprime left $H$-module algebra and let
$\EA$ be the self-injective hull of $A$ as $A^e {\#_\nu} H$-module.
Then the following statements hold:\begin{enumerate}
    \item There exists a bijection between the (essentially) closed
     $H$-stable ideals of $A$,
    the central idempotents of $\EA$ and the central idempotents of $Z(\EA)^H$.
    \item The following statements are equivalent:
    \begin{enumerate}
        \item[(a)] Every direct sum of non-trivial $H$-stable ideals in $A$ is finite.
        \item[(b)] $\EA$ is a finite direct product of $H$-prime $H$-module algebras.
        \item[(c)] $Z(\EA)^H$ is  a finite product of fields.
    \end{enumerate}
    \item $A$ is $H$-prime if and only if $Z(\EA)^H$ is a field.
\end{enumerate}
\end{corollary}

Let $G$ be a group and consider the group ring $H=k[G]$. Let $A$ be an $k$-algebra where $G$ acts on, then
$G$ acts also on $A^e$ and we can form the skew-group ring $A^e \# G$. The $G$-central closure constructed in \cite{wisbauer96}
coincides with our construction of the central closure $\EA$ as self-injective hull of $A^e \#_\nu H$ since,
as mentioned before, $A^e \#_\nu H$ coincides with the ordinary smash product of $A^e$ and $k[G]$ which is
in this case the skew-group ring of $A^e$ and $k[G]$.

Using the results of the last section  
we show that our central closure embeds into the Martindale ring of quotients
of a Hopf-module algebra.
For a left $H$-module algebra $A$ our Martindale ring of quotient $Q_0$ constructed with respect to $\cF_B$ where $B=\MA{A}$ coincides with 
the Martindale ring of quotients constructed by Cohen in \cite{Cohen85}.

\begin{proposition}\label{ZenterOfMartindaleForHopf} Let $H$ be a Hopf algebra over $k$ and let $A$ be a left
$H$-semiprime module algebra with right Martindale ring of quotients $Q_0$. Let $\EA$ be the self-injective hull
of $A$  as $A^e {\#_\nu} H$-module. Assume that one of the following conditions hold:
\begin{enumerate}
    \item[(i)]  $A$ is commutative or
    \item[(ii)] $A$ is semiprime or
    \item[(iii)] $H$ has a bijective antipode
    \end{enumerate}
then
\begin{enumerate}
\item[(1)]
$\psi: Z(\EA)^H \rightarrow Z(Q_0)^H \text{ with } f \mapsto [f : I_f \rightarrow A]$
is an isomorphism of $k$-algebras where $I_f := f^{-1}(A)\cap A$.
\item[(2)] $Z(Q_0)^H$ is a von Neumann regular self-injective $k$-algebra.
\item[(3)] $\overline{\psi}: \EA \rightarrow Q_0 \text{ with } \gamma \mapsto [L_\gamma : I \rightarrow A]$
is an injective homomorphism of $k$-algebras where $I=\bigcup_{j=1}^n I_{x_j}$ for $\gamma = \sum_{j=1}^n a_jx_j$.
\item[(4)] $\EA$ is isomorphic to the subalgebra of $Q_0$ generated by $A$ and $Z(Q_0)^H$.
\end{enumerate}
\end{proposition}

\begin{proof} We just have to check that any of the hypothesis (i-iii) implies that the left and right annihilator
of an $H$-stable ideal in $\cF_H$ of a left $H$-module algebra is $H$-stable.
First of all note that the left annihilator in $A$ of an $H$-stable (left) ideal $I$ is always $H$-stable. Let $a\in A$ satisfy
$aI=0$ then for all $h\in H$ and $x\in I$ one has:
$$ (h\cdot a)x = \sum_{(h)} h_1\cdot \left( a (S(h_2)\cdot x) \right) = 0.$$
We still  have to  show that the right annihilator in $A$ of an essential $H$-stable ideal is $H$-stable.
Case (i): If $A$ is commutative then left and right annihilator are equal and hence $H$-stable.
Case (ii): Let $I$ be an essential $H$-stable ideal. Since the left annihilator $l.ann_A(I)$
of $I$ in $A$ is $H$-stable, $l.ann_A(I)=0$ follows by Lemma \ref{wesentlicheIdeale}. If $A$ is semiprime, then also $r.ann_A(I)=0$ follows,
i.e. the right annihilator of any essential $H$-stable ideal is $H$-stable.\\
Case (iii): If $H$ has a bijective antipode, then the right annihilator in $A$ of an $H$-stable ideal $I$ is always $H$-stable.
 Let $a\in A$ satisfy
$Ia=0$ then for all $h\in H$ and $x\in I$ one has:
$$ x(h\cdot a) = \sum_{(h)} h_2 \cdot \left( (S^{-1}(h_1)\cdot x) a  \right) = 0.$$
Hence we can apply Theorem \ref{CentralClosureEmbedds}and Proposition \ref{ZenterOfMartindale}.
\end{proof}

In \cite{Matczuk92} Matczuk constructs the central closure of an $H$-prime module algebra $A$
directly as the subalgebra of the Martindale
quotient ring $Q_0$ of a module algebra $A$,  generated by $A$ and $Z(Q_0)^H$. We see by $(4)$
 that his construction coincides with ours.

A $H$-semiprime left $H$-module algebra $A$ is called {\bf $H$-centrally closed},
if $\EA=A$ holds. From Theorem \ref{zentraler Abschluss}(4) follows that if $A$ is an $H$-semiprime left $H$-module algebra,
then $\EA$ is $H$-centrally closed.  Is $A$ $H$-prime, then  $\EA$ is $H$-prime and $Z(\EA)^H$ is a field.
We might consider  $\bar{H}:=H \otimes Z(\EA)^H$ as a $Z(\EA)^H$-Hopf algebra and $\EA$ becomes a
$\bar{H}$-prime left $\bar{H}$-module algebra over the field $k:=Z(\EA)^H$.
Was $H$ separable over $k$, then $\bar{H}$ is also separable over $k$  and hence finite dimensional and semisimple.
Thus questions with respect to $H$-prime module algebras over separable Hopf algebras $H$ over commutative rings can be
reduced to $H$-prime $H$-centrally closed module algebras over finite dimensional semisimple Hopf algebras over fields.

\bibliographystyle{amsplain}


\small{Address: Centro de Matem\'{a}tica da Universidade do Porto, Rua Campo Alegre 687, 4169-007 Porto, Portugal; Email: clomp@fc.up.pt}
\end{document}